\documentclass[12pt,draft]{article}
\usepackage{amssymb}
\usepackage{amsmath}
\usepackage{latexsym}

\newcommand{\R}{{\mathbb R}}
\newcommand{\Z}{{\mathbb Z}}

\newcommand{\eps}{{\varepsilon}}
\newcommand{\definition}{\medskip\noindent{\bf Definition\ }}

\newcommand{\estimatetwo}{{\frac{B^2+(m-1)B}{2}}}
\newcommand{\estimatethree}{{\frac{B^3+3(m-1)B^2+2B}{6}}}
\newcommand{\estimatep}{{\frac{(B-1)((B-1)^{p-1}-1)}{2p}+\frac{mB}{2}(p-1)}}

\newtheorem{theorem}{Theorem}
\newtheorem{lemma}{Lemma}[section]
\newtheorem{remark}{Remark}[section]

\newcommand{\beginproof}{\begin{trivlist}
\rm\item[\hspace{\labelsep}{ P\,r\,o\,o\,f. }]}
\newcommand{\proofend}{$\Box$\end{trivlist}}
\author{Fedor Duzhin}
\title{Lower estimates for the number of closed trajectories of generalized billiards}
\date{}
\begin{document}
\maketitle
\section{Introduction}

Closed billiard trajectories is a classical object first considered
by George Birkhoff. A billiard is motion of a particle inside some
domain when field of force is lacking, then the particle moves along
a geodesic line and rebounds from the domain's boundary making the
angle of incidence be equal to the angle of reflection. Closed
trajectories of such a motion are connected with different areas of
mathematics. For example, closed billiard trajectories correspond to
closed geodesics of the following space: we take two copies of the
given domain and glue corresponding points of the boundaries. One can
obtain an another example noticing that the minimal number of closed
billiard trajectories is an invariant of a knot or, say, of a plane
curve.

George Birkhoff stated and solved the following problem in
\cite{birkhoff}: given an integer $k$, estimate from below the
number of closed billiard trajectories with exactly $k$ rebounds.
More precisely, he proved that if $k>2$ is an integer and
$T\subset\R^2$ is a strictly convex domain, then there exist at
least $\varphi(k)$ closed bliiard trajectories with exactly $k$
rebounds. Here $\varphi(k)$ is Euler's function that is an amount of
coprime with $k$ integers not exceeding $k$.

The billiard ball rebounds in this theorem from the boundary of
two-dimensional domain that topologically is a circle embedded in the
Euclidean plane. It's not hard to replace this circle with an
arbitrary manifold embedded in a Euclidean space of any dimension. In
fact, the billiard ball rebounding from something of codimension
greater than $1$ can do it in infinitely many directions, but only
finitely many of them make the ball get back to the manifold.

Many mathematicians tried to estimate the number of closed billiard
trajectories. Morse theory was applied to this problem by Morse
himself. In \cite{morse} Morse investigated the simplest case: a
manifold is an $m$-sphere and closed billiard trajectories consist
of only two points. The best known estimate for the number of closed
billiard trajectories with two rebounds was found by P.~Pushkar in
\cite{pushkar}.

Estimates for the number of closed trajectories of an arbitrary
period for billiards in multi-dimensional convex domains were proved
by M.~Farber and S.~Tabachnikov in \cite{farber} and \cite{farber2}.

In our paper \cite{duzhin} a general estimate for the number of
closed trajectories of period $3$ was obtained. Unfortunately, the
paper \cite{duzhin} contains an error noticed by M.~Farber and
S.~Tabachnikov: in fact, the estimate was proved there only for
manifolds lying in a boundary of some strictly convex domain.

The main goal of the present paper is, first, to correct that
mistake (lemma \ref{morselemma}). Besides, we give a general
estimate for closed trajectories of any prime period, when the
billiard ball rebounds from an arbitrary submanifold of a Euclidean
space (theorem \ref{thp}).

\section{Preliminaries}

Let $M$ be a smooth closed connected $m$-dimensional manifold
embedded in the Euclidean space $\R^n$ (so $m<n$), $k\ge 2$ an
integer. The dihedral group $D_k$ acts on $M$'s Cartesian power
$M^{\times k}$, this action is given by the cyclic permutation
\begin{equation}
(x_1,x_2,\dots,x_k)\rightarrow(x_2,x_3,\dots,x_k,x_1)
\end{equation}
and the reflection
\begin{equation}
(x_1,x_2,\dots,x_k)\rightarrow(x_k,x_{k-1},\dots,x_1).
\end{equation}
The cyclic permutation corresponds to the fact that a closed polygon
may be considered starting from any of its vertices, while the
reflection means that the direction can be reversed. In fact, closed
polygons are points of the quotient space $M^{\times k}/D_k$.

\definition
{An ordered set of points $(x_1,\dots,x_k)\in M^{\times k}$
considered up to the action of the dihedral group $D_k$ is said to
be {\em a closed (periodic, or $k$-periodic}) billiard trajectory if
for any cyclic index $i$ (we mean $i=i+k$) the following conditions
hold:
\begin{enumerate}
\item $\displaystyle x_i\neq x_{i+1}$,
\item $\displaystyle \frac{x_i-x_{i+1}}{\|x_i-x_{i+1}\|}+
\frac{x_i-x_{i-1}}{\|x_i-x_{i-1}\|}\perp T_{x_i}M$.
\end{enumerate}
}

Note that the second of these conditions is the same as the angle of
incidence equals the angle of reflection.

Let us introduce the following notation:
\begin{equation}
\begin{array}{l}
\displaystyle\tilde\Delta=\bigcup_{i\in\Z_k}\{x_i=x_{i+1}\} \subset M^{\times
k},\\
\displaystyle\Delta=\tilde\Delta/D_k
\end{array}
\end{equation}
is the diagonal consisting of all the closed polygons, at least one
of whose segments vanishes. We see that a closed billiard trajectory
is a point of the space $(M^{\times
k}\setminus\tilde\Delta)/D_k=(M^{\times k}/D_k)\setminus\Delta$.

In order to emphasize that $\Delta$ consists of closed polygons with
$k$ segments, we'll write $\Delta_k$. If we need to accentuate that
$\Delta$ consists of closed polygons with vertices belonged to the
manifold $M$, we'll write $\Delta_M$.

Let
\begin{equation}
{l}=\sum_{i\in\Z_k}\|x_i-x_{i+1}\|:M^{\times k}\to\R
\end{equation}
be the length function of a closed polygon, all of whose vertices
lie on the manifold $M$. Obviously, the function $l$ is smooth
outside the diagonal ${\tilde\Delta}$. It's easy to see that $l$ is
invariant under the action of the dihedral group $D_k$, so it
essentially is a function on the quotient space $M^{\times k}/D_k$.

It is well known that closed billiard trajectories with $k$ segments
(or, more exactly, their inverse images under the natural projection
$M^{\times k}\to M^{\times k}/D_k$) are exactly the critical points
of the function $l$ outside of the diagonal $\tilde\Delta$.

\definition A embedding $M\to\R^n$ is {\em generic} (or, more precisely,
{\em $k$-generic}) if all the critical points of all the functions
$l_{k'}$ with $k'<k$ outside of the corresponding diagonals
$\tilde\Delta_{k'}$ are non-degenerate.\footnote{Such embeddings
form an open dense set in the space of all embeddings, for details
see \cite{farber}.}.

Thus our problem is to estimate the minimal number $BT_p(M))$ of
critical points of the function $l$ for all generic embeddings
$M\to\R^n$ (with fixed $M$ and unfixed $n$). We do solve
this problem only for prime $k=p$\footnote{First, if $k=ab$, then
among the all $k$-periodic billiard trajectories there are
$a$-periodic ones repeated $b$ times. Second, the
action of the dihedral group is free only when $k$ is prime, else the
quotient space $(M^{\times k}/D_k)\setminus\Delta$ is not a smooth
manifold.}.

Our main statement is
\begin{theorem}\label{thp}
Let $M$ --- be a smooth closed connected $m$-dimansonal manifold,
$p>3$ is a prime integer\footnote{For $p=2,3$ this estimate can be
strengthened. In fact, the number of closed billiard trajectories of
period $2$ is at least $\estimatetwo$ (see \cite{pushkar}), for
period $3$ the estimate is $\estimatethree$ (see \cite{duzhin}).}.
Put $k_i=\dim H_i(M;\Z_2)$, $B=\sum_{i=0}^m k_i$. Then the minimal
number of closed $p$-periodic billiard trajectories for all generic
embeddings of the manifold $M$ into a Euclidean space satisfy
\begin{equation}
BT_p(M)\ge\estimatep.
\end{equation}
\end{theorem}

\beginproof
By the Morse inequalities (lemma \ref{morselemma}), it follows that
\begin{equation}
BT_p(M)\ge \sum_{i=0}^{pm}\dim H_q(M^{\times p}/D_p,\Delta;\Z_2).
\end{equation}
By the results of the paper \cite{duzhin}, it follows that if the
homology groups of spaces $X_1$ and $X_2$ are isomorphic, then
\begin{equation}
H_*(X_1^{\times p}/D_p,\Delta_{X_1};\Z_2)\cong H_*(X_2^{\times
p}/D_p,\Delta_{X_2};\Z_2)
\end{equation}
as well\footnote{This statement is proved in \cite{duzhin} only for
$p=2,3$, but one can easily generalize it for the case of arbitrary
$p$.}. Hence,
\begin{equation}
\sum_{i=0}^{pm}\dim H_q(M^{\times p}/D_p,\Delta;\Z_2)
=\sum_{i=0}^{pm}\dim H_q(X^{\times p}/D_p,\Delta_X;\Z_2),
\end{equation}
where $X$ is the bouquet of spheres
\begin{equation}
S^m\vee S^{m-1}_{1}\vee \dots\vee S^{m-1}_{k_{m-1}}
\vee\dots\vee S^{1}_{1}\vee\dots\vee S^{1}_{k_{1}}.
\end{equation}
Finally, lemma \ref{bettisumlemma} implies that
\begin{equation}
\sum_{i=0}^{pm}\dim H_q(X^{\times p}/D_p,\Delta_X;\Z_2)\ge
\estimatep.
\end{equation}
This completes the proof. \proofend

\section{Morse inequalities}\label{morsesection}
Let us state the main lemma first and then proceed with all the
propositions needed for its proof.

\begin{lemma}\label{morselemma}
Let $p$ be a prime integer, $M$ be a smooth closed connected
$p$-generic submanifold of the Euclidean space $\R^{n}$. Then there
exist at least
\begin{equation}\label{bettiestimate}
\sum_{q=0}^{mp}\dim H_q(M^{\times p}/D_p,\Delta;\Z_2)
\end{equation}
$p$-periodic billiard trajectories for the
manifold $M$.
\end{lemma}

\beginproof
If $p=2$ or if $M$ lies in a boundary of a strictly convex domain,
then the lemma is proved in \cite{duzhin}.

If $M$ does not, we can deform the embedding $M\to\R^n$ slightly such
that $M$ gets to a boundary of a strictly convex domain. Indeed,
$M\subset\R^{n}\subset\R^{n+1}$ and $\R^n$ can be deformed to a sphere
in $\R^{n+1}$. This deformation is small on $M$ itself. Lemma
\ref{homotopylemma} implies that the number of closed billiard
trajectories remains the same. This completes the proof.\proofend

Since lemma \ref{morselemma} is already proved for $p=2$, in this
section we may assume that $p>2$.

Let us introduce the following functions:
\begin{equation}\label{functionfk}
\begin{array}{l}
\displaystyle f_2=\langle a,x_1-x_2\rangle:S^{n-1}\times M\times M\to\R,\\
\displaystyle f_k=\sum_{i\in\Z_k}\langle
a_i,x_i-x_{i+1}\rangle:(S^{n-1})^{\times k}\times M^{\times k}\to\R,
\ k\ge 3,
\end{array}
\end{equation}
where
$a,a_i\in S^{n-1}=\left\{u_1^2+\dots+u_n^2=1\right\}\subset\R^n$,
$x_i\in M$ and, as above, $i=i+k$. Here $S^{n-1}\subset\R^n$ and
$M\subset\R^n$, that's why all the scalar products are well defined.

\begin{lemma}
Suppose $(x_1,\dots,x_k)$ is a closed billiard trajectory that is a
critical point of the function
\begin{equation}
l_k=\sum_{i\in\Z_k}\|x_i-x_{i+1}\|:M^{\times k}\to\R.
\end{equation}
Then
\begin{equation}
P_0=\left(\frac{x_1-x_{2}}{\|x_1-x_{2}\|},\frac{x_2-x_{3}}{\|x_2-x_{3}\|},\dots,
\frac{x_k-x_{1}}{\|x_k-x_{1}\|},x_1,x_2,\dots,x_k\right)
\end{equation}
is a critical point of the function $f_k$. Similarly,
\begin{equation}
P_0=\left(\frac{x_1-x_{2}}{\|x_1-x_{2}\|},x_1,x_2\right)
\end{equation}
is a critical point of the function $f_2$.
\end{lemma}

\beginproof
Let $P_0=(a_1,\dots,a_k,x_1,\dots,x_k)$ be a critical point of the
function $f_k$ such that $x_i\neq x_{i+1}$ for any
$i\in\Z_k$. We have $\frac{\partial f_k}{\partial
a_i}(P_0)=0$ and
$\frac{\partial f_k}{\partial x_i}(P_0)=0$ for all $i$. The first
condition implies
\begin{equation}
x_i-x_{i+1}\perp T_{a_i}S^{n-1},
\end{equation}
that is $a_i\parallel x_i-x_{i+1}$. The second condition means that
\begin{equation}
a_{i-1}-a_{i}\perp T_{x_i}M.
\end{equation}
Let now $(x_1,\dots,x_k)$ be a critical point of the function
$l_k$. It follows that
\begin{equation}
\frac{x_i-x_{i+1}}{\|x_i-x_{i+1}\|}-
\frac{x_{i-1}-x_i}{\|x_{i-1}-x_i\|}\perp T_{x_i}M.
\end{equation}
Put $a_i=\frac{x_i-x_{i+1}}{\|x_i-x_{i+1}\|}$. This completes the
proof.
\proofend

\begin{lemma}\label{kprimelemma}
Let $P_0=(a_1,\dots,a_k,x_1,\dots,x_k)$ be a critical point of the
function $f_k$ with $k>2$ such that some of $x_i$ coincide. Actually
assume that $\beta_1,\dots,\beta_{k'}$ are integers such that
\begin{enumerate}
\item $\beta_1,\dots,\beta_{k'}\ge 1$,
\item $\beta_1+\dots+\beta_{k'}=k$,
\item $1<k'<k$,
\end{enumerate}
and put
$\alpha_i=\beta_1+\dots+\beta_i$. Suppose $x_i$ to coincide as
follows
\begin{equation}
\begin{array}{lclclcl}
x_1&=&\dots&=&x_{\alpha_1}&\neq\\
x_{\alpha_1+1}&=&\dots&=&x_{\alpha_2}&\neq\\
 &\vdots\\
x_{\alpha_{k'-1}+1}&=&\dots&=&x_{k}&\neq&x_1
\end{array}
\end{equation}
Then the following conditions hold:
\begin{enumerate}
\item the point
$P'_0=(a_1,a_{\alpha_1+1},\dots,a_{\alpha_{k'-1}+1},
x_{\alpha_1},\dots,x_{\alpha_{k'}})$ is critical for the function
$f_{k'}$.
\item the point $P_0$ belongs to a critical manifold
$M_0$ given by
\begin{equation}
\begin{array}{lcllll}
a_2-a_1,&\dots,&a_{\alpha_1}-a_{\alpha_1-1},&a_{\alpha_1+1}-
a_{\alpha_1}&\perp&T_{x_{\alpha_1}}M,\\
a_{\alpha_1+2}-a_{\alpha_1+1},&\dots,&a_{\alpha_2}-
a_{\alpha_2-1},&a_{\alpha_2+1}-a_{\alpha_2}&\perp&T_{x_{\alpha_2}}M,\\
 & \vdots\\
a_{\alpha_{k'-1}+2}-a_{\alpha_{k'-1}+1},&\dots,&a_{\alpha_{k'}}-a_{\alpha_{k'}-1},&
a_{1}-a_{k}&\perp&T_{x_{k}}M.
\end{array}
\end{equation}
\end{enumerate}
\end{lemma}

\beginproof
Without loss of generality consider the simplest case:
$x_1=x_2$ and $x_i\neq x_{i+1}$ for $i\neq 1$. As above, we
obviously have
\begin{equation}
a_i\parallel x_i-x_{i+1},\ \ i\neq 1,
\end{equation}
and
\begin{equation}
a_i-a_{i+1} \perp T_{x_i}M,\ \ i\in\Z_k.
\end{equation}
For $x_1=x_2$ we obtain that $a_k-a_1\perp T_{x_1}M$ and
$a_1-a_2\perp T_{x_1}M$. Summing these two conditions, we get
that $a_k-a_2\perp T_{x_1}M$, while all possible
$a_1$ form an $(n-m-1)$-sphere. For $n=m+1$ this
sphere is just a couple of points, but, in fact, we do not need to
consider this case very detailed, since in further we always have
$n>m+1$.

Notice that for $k'<k-1$ these critical manifold can be products of
spheres.
\proofend

\begin{remark}
Besides, there exists a critical manifold
$M^{(0)}$ given by
\begin{equation}
\begin{array}{l}
x_1=\dots=x_k=x\in M,\\
a_i-a_j\perp T_x M,\ i\neq j.
\end{array}
\end{equation}
It is a bundle over $M$, the fiber is defined by the second of these
conditions.
\end{remark}

\begin{lemma}
Let $(x_1,\dots,x_k)$ be a non-degenerate critical point of the
function $l_k$ and $\mu$ be its Morse index. Then the corresponding
critical point
\begin{equation}
P_0=(a_1,\dots,a_k,x_1,\dots,x_k),\
a_i=\frac{x_i-x_{i+1}}{\|x_i-x_{i+1}\|},
\end{equation}
of the function $f_k$ is also non-degenerate and its Morse index
equals
$\mu+k(n-1)$ for $k>2$ or $\mu+n-1$ for $k=2$.
\end{lemma}

\beginproof
Assume that $P=(b_1,\dots,b_k,y_1\dots,y_k)$ lies in a small
neighborhood of the critical point $P_0$ being considered. Let us
introduce coordinates in this neighborhood in the following way.
Suppose that
\begin{equation}
y_i=y_i(t_i),
\end{equation}
where $t_i\in\R^m$ is some parametrization for
$y_i$. Put
\begin{equation}
b_i=b_i(s_i),
\end{equation}
where the parametrization $s_i$ for $b_i$ is defined as follows. Let
$A_i$ be an orthogonal operator
$\R^{n-1}\to(y_i-y_{i+1})^\perp$ and $s_i\in\R^{n-1}$ be our
parameter. Put
\begin{equation}
b_i=\frac{\|s_i\|A_is_i+\frac{y_i-y_{i+1}}{\|y_i-y_{i+1}\|}}
{\left\|\|s_i\|A_is_i+\frac{y_i-y_{i+1}}{\|y_i-y_{i+1}\|}\right\|}=
\frac{\|s_i\|A_is_i+\frac{y_i-y_{i+1}}{\|y_i-y_{i+1}\|}}{\sqrt{1+\|s_i\|^2}}
.
\end{equation}
Evidently, $\|b_i\|=1$ and
$b(0)=\frac{y_i-y_{i+1}}{\|y_i-y_{i+1}\|}$. Since $A_is_i\perp
y_i-y_{i+1}$, we have
\begin{equation}
\langle b_i,y_i-y_{i+1}
\rangle=\frac{\|y_i-y_{i+1}\|}{\sqrt{1+\|s_i\|^2}}.
\end{equation}
Thus the following condition holds:
\begin{equation}
\begin{array}{l}
\displaystyle f_k(P)=\sum_{i\in\Z_k}\langle b_i,y_i-y_{i+1}
\rangle=\sum_{i\in\Z_k}\frac{\|y_i-y_{i+1}\|}{\sqrt{1+\|s_i\|^2}}=\\
\displaystyle=l_k(y_1,\dots,y_k)-\frac{1}{2}\sum\|s_i\|^2+\dots
\end{array}
\end{equation}
This concludes the proof.\proofend

\begin{remark}
Suppose $x_i=x_{i+1}$ for some $i$. As we have showed above, in this
case there are some critical manifolds corresponding to critical
points of the function $f_{k'}$, $k'<k$. If
$(x_1,\dots,x_{k'})$ is a non-degenerate critical point of a
function $l_{k'}$, then the corresponding critical manifold is also
non-degenerate. The critical manifold $M^{(0)}$ defined by
$x_1=x_2=\dots=x_k$ is non-degenerate as well.
\end{remark}

\begin{lemma}\label{homotopylemma}
Suppose $F:M\times[0,1]\to\R^n$ is a smooth homotopy such that for
every $t$ the embedding $F_t:M\to\R^n$ is generic, $k\ge 2$ an
integer. Then the homotopy $F$ keeps the number of closed billiard
trajectories with $k$ vertices.
\end{lemma}

\beginproof
Consider the homotopy $F:M\times[0,1]\to\R^n$. Denote $f_k$- and
$l_k$-functions corresponding to an embedding $F_t$ by
$f_{kt}$ and $l_{kt}$.

From the previous statements we know that the whole picture is as
follows. Closed billiard trajectories (those are non-degenerate
isolated critical points of the function $l_k$) correspond to
non-degenerate critical points of the function $f_k$ that is a smooth
function defined on the smooth manifold
\begin{equation}
(S^{n-1})^{\times k}\times M^{\times k}.
\end{equation}
We suppose every embedding $F_t:M\to\R^n$ to be generic, thus
$f_{kt}$ has an amount of isolated critical points and several
non-degenerate critical manifolds corresponding to isolated critical
points of the functions $f_{k't}$ with $k'<k$.

Thus when the embedding $F_0:M\to\R^n$ is being deformed, isolated
critical points of the function $f_{k}$ could disappear and be born
only from non-degenerate critical manifolds that is impossible.
Indeed, suppose an isolated critical point is born at
$t=t_0$. We mean that there exists $M_t$ --- a non-degenerate
critical manifold of the function $f_{kt}$ for $|t-t_0|$ small enough
and for $t>t_0$ there exists an isolated critical point
$P_t$ such that $\lim_{t\to t_0+0}P_t=P_0\in M_{t_0}$. By
Morse-Bott theory, there are coordinates $t,X^1,\dots,X^N$ in a
neighborhood
$U\subset (S^{n-1})^{\times k}\times
M^{\times k}\times[t_0-\eps,t_0+\eps]$  of the point $P_{t_0}$ such
that
$M_t$ is given by
$X^1=\dots=X^r=0$ and
\begin{equation}
f_{kt}=C(t)-(X^{r+1})^2-\dots-(X^{r+s})^2+(X^{r+s+1})^2+\dots+(X^{N})^2.
\end{equation}
We see that in the neighborhood $U$ there are no other isolated
critical points of the function $f_{kt}$. This contradiction
completes the proof. \proofend

\section{Computations for a bouquet of spheres}
\begin{lemma}\label{poincarelemma}
Let $M$ be a smooth closed connected $m$-dimensional manifold,
$k_i=\dim H_i(M;\Z_2)$, $i=0,1,\dots,m$, $B=\sum_{i=0}^m k_i.$ Then
\begin{equation}
\sum_{i=1}^mik_i=\frac{mB}{2}.
\end{equation}
\end{lemma}
\newpage
\beginproof
Poincar\'e duality implies that:
\begin{equation}
\begin{array}{l}\displaystyle
\sum_{i=1}^m ik_i=
\sum_{i=0}^m ik_i=
\frac{1}{2}\sum_{i=0}^m\left(ik_i+\left(m-i\right)k_{m-i}\right)=\\
\displaystyle\frac{1}{2}\sum_{i=0}^m\left(ik_i+\left(m-i\right)k_{i}\right)=
\frac{mB}{2}.
\end{array}
\end{equation}
\proofend

\begin{lemma}\label{bettisumlemma}
Let $M$ be a smooth closed connected $m$-dimensional manifold, p and
odd prime, $k_i=\dim H_i(M;\Z_2)$,
$B=\sum_{i=0}^m k_i$.
Suppose
\begin{equation}
 X=S^m\vee S^{m-1}_{1}\vee \dots\vee S^{m-1}_{k_{m-1}}
\vee\dots\vee S^{1}_{1}\vee\dots\vee S^{1}_{k_{1}}.
\end{equation}
Then
\begin{equation}\label{bettisum}
\sum_{i=1}^{pm}\dim H_i(X^{\times p}/D_p,\Delta_X;\Z_2)\ge
\frac{(B-1)((B-1)^{p-1}-1)}{2p}+\frac{mB}{2}(p-1).
\end{equation}
\end{lemma}

\beginproof
Consider the bouquet of spheres $X$. By $X_0$ denote the common point
of all the spheres. Let $X_i$ be the $i$th sphere of the bouquet
without the point $X_0$, so topologically $X_i$ is a Euclidean space
$\R^q$ and $X=X_0\cup X_1\cup\dots\cup X_{B-1}$ is a cell decomposition.

Clearly we have
\begin{equation}
X^{\times p}=\bigcup_{i_1,\dots,i_p}X_{i_1}\times\dots\times X_{i_p}
\end{equation}
is a cell decomposition of the Cartesian power $X^{\times p}$. What
we do need is to construct its subdecomposition such that
\begin{itemize}
\item it is invariant under the action of the dihedral group $D_p$,
\item the diagonal $\Delta$ is a cell subspace.
\end{itemize}
Note that if $i_1\neq i_2\neq i_3\neq\dots\neq i_p\neq i_1$, then
$X_{i_1\dots i_p}=X_{i_1}\times\dots\times X_{i_p}$ does not intersect
the diagonal. It follows that $X_{i_1\dots i_p}$ is a cell of the
decomposition being constructed and its boundary is zero.

Consider now $X_{i_1\dots i_p}=X_{i_1}\times\dots\times X_{i_p}$ such
that $i_\alpha=i_{\alpha+1}$ for some $\alpha\in\Z_p$. First suppose
that not all of the $i_\alpha$ coincide. Without loss of generality
we can assume that
$i_1=\dots=i_{\beta_1}\neq i_{\beta_1+1}=\dots=i_{\beta_2}
\neq\dots\neq i_{\beta_u+1}=\dots=i_{p}\neq i_1$.
We construct a cell subdecomposition for all
$X_i\times X_i\times\dots\times X_i$
and the decomposition for the whole
$X_{i_1\dots i_p}$ would be their tensor product.

Each $X_i$ is topologically a Euclidean space $\R^q$. Thus we deal
with the Cartesian power $(\R^q)^{\times\beta}$. Let the $j$th
$\R^q$ have coordinates $x_1^j,\dots,x_q^j$. A cell is given by
the following conditions:
\begin{equation}
\begin{array}{lllllll}
x^1_1&\eps^1_1&x^2_1&\eps^2_1&\dots&\eps^{\beta-1}_1&x^\beta_1,\\
x^1_2&\eps^1_2&x^2_2&\eps^2_2&\dots&\eps^{\beta-1}_2&x^\beta_2,\\
\vdots\\
x^1_{q}&\eps^1_{q}&x^2_{q}
&\eps^2_{q}&\dots&\eps^{\beta-1}_{q}&x^\beta_{q}.
\end{array},
\end{equation}
where each $\eps^*_*$ is one of the signs $<$, $>$, or $=$.

Now consider $X_{i_1\dots i_p}$ having $i_1=\dots=i_p=i$. Then the
cell subdecomposition for this thing is given by the same
construction with inequalities $x^{\beta}_j\eps^{\beta}_j x^1_j$
added. Clearly, $\eps^{1}_j,\eps^{2}_j,\dots,\eps^{\beta}_j$ should
not be all $<$ or all $>$, since in this case the system of
inequalities has no solutions at all.

We have just constructed the cell decomposition for the space
$X^{\times p}$. Denote the corresponding chain complex by
$C(X^{\times p})$. It induces the cell decomposition for the quotient
$X^{\times p}/D_p$ with the diagonal contracted to a point. Let us
denote the induced chain complex by $C(X^{\times p}/D_p,\Delta_X)$.
Our goal is to calculate its homology
\begin{equation}
H_*C(X^{\times p}/D_p,\Delta_X).
\end{equation}

First consider $X_{i\dots i}$ for some fixed $i>0$. Suppose $\dim
X_i=q$. Note that all $X_{i_1\dots i_p}$ such that
$i_\alpha$ is either $0$ or $i$ for all $\alpha=1,\dots,p$
form a chain subcomplex. Denote it by $C^+(X_{i\dots i})$. Moreover,
there is no cell outside $C^+(X_{i\dots i})$ such that its algebraic
boundary contains terms lying in $C^+(X_{i\dots i})$. Hence
$C^+(X_{i\dots i})$ is a direct summand
in $C(X^{\times p}/D_p,\Delta_X)$. Obviously, $C^+(X_{i\dots i})$
coincides with a chain complex for a sphere $C((S^q)^{\times
p}/D_p,\Delta_{S^q})$. By the results of M.~Farber and S.~Tabachnikov
(see \cite{farber}, \cite{farber2}), it follows that
\begin{equation}
\sum\dim H_\alpha(C^+(X_{i\dots i});\Z_2)=q(p-1).
\end{equation}
Summing for all $i$ and using lemma \ref{poincarelemma}, we obtain
that the contribution to the sum \ref{bettisum} being calculated
equals
\begin{equation}
\frac{mB}{2}(p-1).
\end{equation}

Now consider $X_{i_1\dots i_p}=X_{i_1}\times\dots\times X_{i_p}$ for
$i_1\neq i_2\neq i_3\neq\dots\neq i_p\neq i_1$.
Each of these $X_{i_1\dots i_p}$ is a cell such that
\begin{itemize}
\item its algebraic boundary is zero,
\item it is not contained in an algebraic boundary of any other
cell.
\end{itemize}
Hence it forms a chain subcomplex in $C(X^{\times p}/D_p,\Delta_X)$
consisting of only one group with only one generator and zero
boundary operator. Let us denote this chain complex by $C(X_{i_1\dots
i_p})$. It contributes $1$ to the sum \ref{bettisum}.

It is well known from combinatorics that the number of all
$X_{i_1\dots i_p}$ having
$i_1\neq i_2\neq i_3\neq\dots\neq i_p\neq i_1$ equals
\begin{equation}
(B-1)((B-1)^{p-1}-1).
\end{equation}
Anyway let us prove it. Suppose $N(p)$ is the number of all
$p$-tuples $(i_1,\dots,i_p)$ such that $0\le i_\alpha\le B-1$ and
$i_1\neq i_2\neq i_3\neq\dots\neq i_p\neq i_1$. Let now $p$ be not
necessarily prime. Then we have
\begin{equation}
N(p)=B(B-1)^{p-1}-N(p-1).
\end{equation}

Indeed, $i_1$ may be chosen in $B$ ways. Each $i_\alpha$,
$\alpha=2,\dots,p$ may be chosen in $B-1$ ways to be different from
$i_{\alpha-1}$. It gives $B(B-1)^{p-1}$. If $i_1=i_p$,
then
$(i_1,i_2,\dots,i_{p-1})$ is a correct $(p-1)$-tuple.
Thus recalling $p$ is an odd prime we get
\begin{equation}
\begin{array}{l}
\displaystyle N(p)=B(B-1)^{p-1}-N(p-1)=\\
\displaystyle B(B-1)^{p-1}-B(B-1)^{p-2}+N(p-2)=\\
\displaystyle B(B-1)^{p-1}-B(B-1)^{p-2}+B(B-1)^{p-3}-\dots-B(B-1)=\\
\displaystyle B(1-B)\frac{1-(1-B)^{p-1}}{1-(1-B)}=(B-1)((B-1)^{p-1}-1).
\end{array}
\end{equation}
This calculation gives us
\begin{equation}
\frac{(B-1)((B-1)^{p-1}-1)}{2p}
\end{equation}
after factorizing by the action of the dihedral group $D_p$ .
\proofend

\end{document}